\documentclass[12pt]{article}
\usepackage{amsmath,amsfonts,amssymb,theorem}

{\setlength{\labelwidth}{3.5cm}\setlength{\labelsep}{1.00em}
\setlength{\leftmargin}{5.00cm}\setlength{\rightmargin}{3.00cm}

{\setlength{\labelwidth}{3.5cm}\setlength{\labelsep}{1.5em}
\setlength{\leftmargin}{5.00cm}\setlength{\rightmargin}{3.00cm}

\addtolength{\hoffset}{0.25cm}

 \usepackage[matrix, arrow]{xy}

\DeclareMathOperator{\LCS}{LCS}

\DeclareMathOperator{\lct}{lct}

 \numberwithin{equation}{subsubsection}
 \numberwithin{footnote}{section}

 \newtheorem{cor}[subsection]{Corollary}

 \newtheorem{thm}[subsection]{Theorem}
 \newtheorem{conj}[subsection]{Conjecture}

{
 \theorembodyfont{\rmfamily}
 \newtheorem{defn}[subsection]{Definition}
 
 \newtheorem{exa-cr}[subsubsection]{Example--Construction}

 \newtheorem{rem}[subsection]{Remark}

}
 \newcommand{\qed}{\ifhmode\unskip\nobreak\fi\quad\ensuremath\square}
 \newenvironment{proof}{\paragraph{Proof}}{\par\medskip}

 \newcommand{\ke}[1]{$\acute{\mbox{e}}$}
 \newcommand{\ku}[1]{$\acute{\mbox{u}$}}
 \newcommand{\kl}[1]{$\acute{\mbox{l}}$}
 \newcommand{\kh}[1]{$\acute{\mbox{h}}$}
 \newcommand{\kr}[1]{$\acute{\mbox{r}}$}
 \newcommand{\kx}[1]{$\acute{\mbox{x}}$}
 \newcommand{\ki}[1]{${\^\i}$}

\title{ \normalsize{\textbf{
ACC for log canonical thresholds and termination of log flips
}}}

\small

\author{\textbf{ Caucher Birkar\thanks{Supported by the EPSRC.}}}

\begin{document}

\maketitle

\begin{abstract} 
We prove that the ascending chain condition (ACC) for log canonical (lc) thresholds in dimension $d$ and Special Termination in dimension $d$ imply the termination of any sequence of log flips starting with a $d$-dimensional lc pair of nonnegative Kodaira dimension. In particular, in characteristic zero, the latter termination in dimension 4 follows from Alexeev-Borisov's conjecture in dimension 3. 
\end{abstract}

\tableofcontents


\section{Introduction}

Following the fundamental work of Shokurov [Sh3] on the existence of log flips,  one of the main open problems in log minimal model program (LMMP) is the termination of log flips, in particular, in dimension 4. 

\begin{conj}[Termination of log flips]\label{termination}
Any sequence of log flips/$Z$, with respect to a lc log divisor $K_X+B$, terminates. 
\end{conj}

In this paper we prove the following

 \begin{thm}\label{theorem}
 ACC for lc thresholds in dimension $d$ and Special Termination in dimension $d$ imply the termination of any sequence of log flips/$Z$ starting with a $d$-dimensional lc pair $(X_1,B_1)$ of nonnegative Kodaira dimension. 
\end{thm}

\begin{cor}\label{corollary}
In characteristic zero, Alexeev-Borisov's Conjecture in dimension $3$ implies the termination of any sequence of log flips/$Z$ starting with a $4$-dimensional lc pair $(X_1,B_1)$ of nonnegative Kodaira dimension.
\end{cor}

Here we recall the  ACC for lc thresholds Conjecture and Special Termination Conjecture due to Shokurov and the Alexeev-Borisov's Conjecture due to Alexeev, A. Borisov and L. Borisov. For a set $S\subseteq \mathbb{R}$, $B\in S$ means that all coefficients of $B$ are in $S$.

\begin{conj}[ACC for lc thresholds]
Suppose that $\Gamma\subseteq [0,1]$ satisfies the descending chain condition (DCC) and  $S\subseteq \mathbb{R}$ is finite. Then the set 
$$\{\lct(M,X,B)|\mbox{ $(X,B)$ is lc of dimension $d$, $B\in \Gamma$ and $M\in S$}\}$$

{\flushleft satisfies} the  ACC where $M$ is an $\mathbb{R}$-Cartier divisor on $X$ and $\lct(M,X,B)$ is the lc threshold of $M$ with respect to $(X,B)$.
\end{conj}

\begin{conj}[Special Termination]
Let $X_i\dasharrow X_{i+1}/Z$ be a sequence of flips/$Z$ starting from a lc pair $(X_1,B_1)$ of dimension $d$. Then there is $I\in \mathbb{N}$ such that the flipping locus does not intersect the locus of lc singularities $\LCS(X_i,B_i)$ for any $i\geq I$. 
\end{conj}

\begin{conj}[Alexeev-Borisov's]
Let $\delta >0$ be a real number. Then,  varieties $X$ for which $(X/\mathrm pt.,B)$ is a d-dimensional $\delta$-lc weak log Fano (WLF) pair for a boundary $B$ form a bounded family.
\end{conj}

The Termination Conjecture (\ref{termination}) was proved in the 3-dimensional terminal case by Shokurov [Sh4], the 3-dimensional Klt and the  4-dimensional terminal cases by Kawamata [K, KMM] and the  4-dimensional canonical case by Fujino [F]. Alexeev also has proved some special cases [A]. But the main development toward a general solution  is the paper by Shokurov [Sh1] which proves that ACC and lower-semicontinuity for minimal log discrepancies (mld's) in dimension $d$ imply the Termination Conjecture in dimension $d$. In particular, ACC for mld's in dimension 4 is enough to prove the Termination in dimension 4 [Sh1].  However, dealing with mld's is not an easy task. Alternatively, in this paper, we use the ACC for lc  thresholds to prove the Termination in the case of nonnegative Kodaira dimension. Moreover, our results also hold in the nonprojective case i.e. the varieties in the sequence are not necessarily globally projective. Recall that pairs with nonnegative Kodaira dimension are exactly those which the LMMP predicts to have a \emph{log minimal model}.

The proof of our results are short and simple. Partly because we use important results of Shokurov, M$^c$Kernan and Prokhorov, and partly because our approach  may be a more natural one. 

\begin{rem}\label{remark}
In characteristic zero, Special Termination in dimension $d$  follows from the LMMP in dimension $\leq d-1$ and the existence of log flips in dimension $d$ [Sh1]. In particular, Special Termination in dimension $4$ is proved [Sh1, cor 4]. Using Special Termination and the arguments of M$^c$Kernan and Prokhorov [MP] one can prove that the ACC for lc thresholds in dimension $4$ follows from the Alexeev-Borisov's Conjecture in dimension $3$. 
\end{rem}

\section{Notations and conventions}

We assume the base field $k$ to be an algebraically closed field. A \emph{log pair} $(X,B)$ consists of a normal variety $X$ over $k$ and $B$ an $\mathbb{R}$-Cartier divisor on $X$ with coefficients in $[0,1]$. More generally, $(X/Z,B)$ consists of a log pair $(X,B)$ equipped with a contraction $X\to Z$ onto a normal variety $Z$ over $k$. We use the usual definitions of singularities of pairs such as Klt, dlt and lc [Sh2], and $\LCS(X,B)$ denotes the locus of lc singularities of $(X,B)$.  

\begin{defn}
Let $(X/Z,B)$ be a pair. We say that the \emph{Kodaira dimension} of $(X/Z,B)$ is nonnegative if there is an effective $\mathbb{R}$-divisor $M$ such that $K_X+B\sim_{\mathbb{R}}M/Z$.
\end{defn}

 Note that if $K_X+B$ is a $\mathbb{Q}$-divisor and the \emph{usual} Kodaira dimension of $(X,B)$ is nonnegative (i.e. $|m(K_X+B)|\neq \emptyset$ for some $m\in\mathbb{N}$), then the Kodaira dimension of $(X/Z,B)$ is also nonnegative  by our definition.

\begin{defn}
Let $M$ be an $\mathbb{R}$-Cartier divisor on $X$. The \emph{lc threshold} of $M$ with respect to a lc pair $(X,B)$ is the real number $\lct(M,X,B):=t$ such that $(X,B+tM)$ is lc but not Klt. 
\end{defn}

A \emph{sequence of flips with respect to $K_{X_1}+B_1$ } is a sequence  $X_i\dasharrow X_{i+1}/Z$ of $K_{X_i}+B_i$-flips$/Z$ where $K_{X_i}+B_i$ is the birational transform of $K_{X_1}+B_1$.

\section{Proof of the main theorem}

\begin{proof}(of Theorem \ref{theorem})

Let  $X_i\dasharrow X_{i+1}/Z$ be a sequence of flips with respect to $K_{X_1}+B_1$. Since the   Kodaira dimension of $(X/Z,B)$ is nonnegative, by definition  $K_{X_1}+B_1\sim_{\mathbb{R}}M_1/Z$ for some effective $\mathbb{R}$-Cartier divisor $M_1$.\\

Step 1.  Let $t_1$ be the lc threshold of $M_1$ with respect to  $K_{X_1}+B_1$. By Special Termination, there is $I\in \mathbb{N}$ such that $\LCS(X_i,B_i+t_1M_i)$ does not intersect the flipping locus for any $i\geq I$ where $M_i$ is the birational transform of $M_1$.\\

Step 2. Now replace  $(X_i,B_i)$ with $(X_{i},B_{i})-\LCS(X_{i},B_{i}+t_1M_i)$ for $i\geq I$. So, we obtain a new sequence of flips. Moreover, the new pair $(X_{I},B_{I}+t_1M_I)$ is Klt. Note that, either the sequence of flips stops at $X_I$ or $M_I\neq 0$ on $X_I$.\\

Step 3. Let $t_2$ be the lc thresholds of $M_I$ with respect to $K_{X_I}+B_I$. Since  $(X_{I},B_{I}+t_1M_I)$ is Klt, we deduce that $t_2>t_1$. Similar to step 1, Special Termination implies that there is $I_1$ such that $\LCS(X_i,B_i+t_2M_i)$ does not intersect the flipping locus for any $i\geq I_1$.\\  

Step 4. By repeating the above process, either the sequence of flips stops or we get an increasing sequence $t_1<t_2<t_3<\dots$ of lc thresholds. The latter contradicts the ACC for lc thresholds.\\ 

$\Box$

\end{proof}

\begin{proof}(of Corollary \ref{corollary})
Immediate by Theorem \ref{theorem} and Remark \ref{remark}.
$\Box$
\end{proof}

\begin{rem}
As mentioned before, our method also works in the nonprojective case ( in particular, for 4-folds) if the same termination holds in the nonprojective case in lower dimensions. But it is not clear yet how to deal with the negative Kodaira dimension case. 
\end{rem}

\begin{rem}
Our method gives a new proof of the termination in dimension $3$ in the nonnegative Kodaira dimension case. Note that our proof, unlike other proofs [K, Sh2], does not use the classification of $3$-fold terminal singularities. 
\end{rem}

\begin{rem}
A plan to attack the Alexeev-Borisov's Conjecture in dimension $3$, using the theory of complements,  is exposed in [B].
\end{rem}

\section{Acknowledgment}

I thank Prof Shokurov for his comments and discussions.

\flushleft
Caucher Birkar\\
Mathematics Institute,\\
Warwick University,\\
Coventry, CV4 7AL\\
UK\\

E-mail: birkar@maths.warwick.ac.uk

\end{document}